\documentclass[preprint,review,10pt]{elsarticle}
\usepackage{amsfonts}
\usepackage{amsmath,amssymb}
\usepackage{graphicx}
\usepackage{setspace}
\usepackage[left=2cm,top=2cm,right=2cm]{geometry}

\newtheorem{definition}{Definition}

\newtheorem{lemma}{Lemma}

\newtheorem{remark}{Remark}

\newtheorem{theorem}{Theorem}
\numberwithin{equation}{section}
 \journal{\textsc{International Journal of Analysis}(Accepted)}
 
\begin{document}
\title{On Simultaneous Approximation of Modified Baskakov Durrmeyer Operators}
\author[label*,label a,label b]{Prashantkumar Patel}
\ead{prashant225@gmail.com}
\author[label a,label c]{Vishnu Narayan Mishra}
\ead{vishnu\_narayanmishra@yahoo.co.in;
vishnunarayanmishra@gmail.com}
\address[label a]{Department of Applied Mathematics and Humanities,
Sardar Vallabhbhai National Institute of Technology, Ichchhanath
Mahadev Dumas Road, Surat-395 007 (Gujarat), India}
\address[label b]{Department of Mathematics,
St. Xavier College, Ahmedabad-380 009 (Gujarat), India}
\address[label c]{L. 1627 Awadh Puri Colony Beniganj, Phase-III, Opposite-Industrial Training Institute (I.T.I.), Ayodhya Main Road,
Faizabad, Uttar Pradesh 224 001, India}
\fntext[label*]{Corresponding authors}
\begin{abstract}
In this manuscript, we  discuss properties of modified
Baskakov-Durrmeyer-Stancu (BDS) operators with parameter
$\gamma>0$. We compute the moments of these modified operators.
Also, establish point-wise convergence, Voronovskaja type
asymptotic formula and an error estimation in terms of second
order modification of continuity of the function for the operators
$B_{n,\gamma}^{\alpha,\beta}(f,x)$.
\end{abstract}

\maketitle
\section{Introduction}
For $x\in [0,\infty)$, $\gamma>0$, $0\le \alpha \le \beta $ and
$f\in C[0,\infty)$, we consider a certain integral type
generalized Baskakov operators as
\begin{eqnarray}\label{5.Eq1}
B_{n,\gamma}^{\alpha,\beta}(f(t),x)&=& \sum_{k=1}^{\infty} p_{n,k,\gamma}(x) \int_0^{\infty}b_{n,k,\gamma}(t) f\left(\frac{nt+\alpha}{n+\beta}\right)dt+ p_{n,0,\gamma}(x)f\left(\frac{\alpha}{n+\beta}\right)\\
&=& \int_0^{\infty} W_{n,\gamma}(x,t)f\left(\frac{nt+\alpha}{n+\beta}\right)dt,\nonumber
\end{eqnarray}
where $$p_{n,k,\gamma}(x)= \frac{\Gamma(n/\gamma + k)}{\Gamma(k+1)\Gamma(n/\gamma)}\cdot\frac{(\gamma x)^k}{(1+\gamma x)^{(n/\gamma)+k}},$$
$$b_{n,k,\gamma}(t)= \frac{\gamma\Gamma(n/\gamma + k+1)}{\Gamma(k)\Gamma(n/\gamma+1)}\cdot\frac{(\gamma t)^{k-1}}{(1+\gamma t)^{(n/\gamma)+k+1}}$$
and
$$ W_{n,\gamma}(x,t)= \sum_{k=1}^{\infty} p_{n,k,\gamma}(x)b_{n,k,\gamma}(t)+ (1+\gamma x)^{-n/\gamma}\delta(t),$$
$\delta(t)$ being the Dirac delta function.\\
\indent The operators defined by (\ref{5.Eq1}) are the
generalization of the integral modification of well-known Baskakov
operators having weight function of some beta basis function. As a
special case, i.e. $\gamma=1$, the operators (\ref{5.Eq1})
 reduce to the operators very recently studied in
 \cite{gupta2012simultaneous,verma2012some}. Inverse results of same type of operators were establish in
 \cite{mishra2013inverse}.
  Also, if $\alpha=\beta= 0$ , the operators (\ref{5.Eq1}) reduce to the operators recently studied in \cite{Vijaygupta2009}
  and if $\alpha=\beta=0$ and $\gamma=1$, the operators (\ref{5.Eq1}) reduce to the operators studied in
  \cite{GuptaNoorBenwal2006}. The $q$-analog of the operators \eqref{5.Eq1} is discussed in \cite{mishra2013approximation}.
 We refer to some of the
important papers on the recent development on similar type of the
operators
\cite{mishra2013hypergeometric,mishra2012simultaneous,patel2015approximation}.
The present paper deals with the study of simultaneous
approximation for the operators $B_{n,\gamma}^{\alpha,\beta}$.
\section{Moments and recurrence relations}
\begin{lemma}\label{5.lemma4}
If we define the central moments, for every $m \in \mathbb{N}$ as
\begin{eqnarray*}
\mu_{n,m,\gamma}(x) = B_{n,\gamma}^{\alpha,\beta}((t-x)^m,x)= \sum_{k=1}^{\infty} p_{n,k,\gamma}(x)
            \int_0^{\infty}b_{n,k,\gamma}(t) \left(\frac{nt+\alpha}{n+\beta}-x\right)^m dt+ p_{n,0,\gamma}(x)\left(\frac{\alpha}{n+\beta}-x\right)^m,
\end{eqnarray*}
then $\mu_{n,0,\gamma}(x)=1$,  $\mu_{n,1,\gamma}(x)= \frac{\alpha-\beta x}{n+\beta}$ and for $n>\gamma m$ we have the following recurrence relation:\\
$\displaystyle
(n-\gamma m)(n+\beta)\mu_{n,m+1,\gamma}(x)= n x (1+\gamma x) \left\{\mu_{n,m,\gamma}^{(1)}(x)+m\mu_{n,m-1,\gamma}(x)\right\}\\
\indent \indent \indent\indent \indent \indent \indent \indent \indent+ \left\{mn+ n^2x- \left(2\gamma m-n\right)\left(\alpha -(n+\beta)x\right)\right\}\mu_{n,m,\gamma}(x)\\
\indent \indent \indent \indent \indent \indent\indent \indent \indent + \left\{  m \gamma (n+\beta) \left(\frac{\alpha}{n+\beta}-x \right)^2-m n \left(\frac{\alpha}{n+\beta}-x\right) \right\}\mu_{n,m-1,\gamma}(x).
$\\
From the recurrence relation, it can be easily verified that for
all $x\in [0,\infty)$, we have $\displaystyle \mu_{n,m,\gamma}(x)=
O(n^{-[(m+1)/2]})$, where $[\alpha]$ denotes the integral part of
$\alpha.$
\end{lemma}
\textbf{Proof:} Taking derivative of above
\begin{eqnarray*}
\mu_{n,m,\gamma}^{(1)}(x)&=& - m \sum_{k=1}^{\infty} p_{n,k,\gamma}(x) \int_0^{\infty}b_{n,k,\gamma}(t) \left(\frac{nt+\alpha}{n+\beta}-x\right)^{m-1} dt -mp_{n,0,\gamma}(x)\left(\frac{\alpha}{n+\beta}-x\right)^{m-1}\\
&&+ \sum_{k=1}^{\infty} p_{n,k,\gamma}^{(1)}(x) \int_0^{\infty}b_{n,k,\gamma}(t) \left(\frac{nt+\alpha}{n+\beta}-x\right)^m dt+ p_{n,0,\gamma}^{(1)}(x)\left(\frac{\alpha}{n+\beta}-x\right)^m\\
&=& -m\mu_{n,m-1,\gamma}(x)+ \sum_{k=1}^{\infty} p_{n,k,\gamma}^{(1)}(x) \int_0^{\infty}b_{n,k,\gamma}(t) \left(\frac{nt+\alpha}{n+\beta}-x\right)^m dt+ p_{n,0,\gamma}^{(1)}(x)\left(\frac{\alpha}{n+\beta}-x\right)^m,
\end{eqnarray*}
$\displaystyle
x(1+\gamma x) \left\{\mu_{n,m,\gamma}^{(1)}(x)+ m\mu_{n,m-1,\gamma}(x)\right\} = \sum_{k=1}^{\infty}x(1+\gamma x) p_{n,k,\gamma}^{(1)}(x) \int_0^{\infty}b_{n,k,\gamma}(t) \left(\frac{nt+\alpha}{n+\beta}-x\right)^m dt\\
\indent \indent \indent\indent \indent \indent\indent \indent \indent\indent \indent \indent+ x(1+\gamma x) p_{n,0,\gamma}^{(1)}(x)\left(\frac{\alpha}{n+\beta}-x\right)^m.
$\\
Using $x(1+\gamma x)p_{n,k,\gamma}^{(1)}(x)= (k-nx)p_{n,k,\gamma}(x)$, we get
\begin{eqnarray}\label{5.eq2.1}
&&x(1+\gamma x) \left\{\mu_{n,m,\gamma}^{(1)}(x)+ m\mu_{n,m-1,\gamma}(x)\right\}\nonumber\\
&&\indent \indent =  \sum_{k=1}^{\infty}(k-n x) p_{n,k,\gamma}(x) \int_0^{\infty}b_{n,k,\gamma}(t) \left(\frac{nt+\alpha}{n+\beta}-x\right)^m dt+ (-nx) p_{n,0,\gamma}(x)\left(\frac{\alpha}{n+\beta}-x\right)^m\nonumber\\
&&\indent \indent = \sum_{k=1}^{\infty}k p_{n,k,\gamma}(x) \int_0^{\infty}b_{n,k,\gamma}(t) \left(\frac{nt+\alpha}{n+\beta}-x\right)^m dt-nx\mu_{n,m,\gamma}(x):= I-nx\mu_{n,m,\gamma}(x).
\end{eqnarray}
We can write $I$ as
\begin{eqnarray}
I&=& \left\{\sum_{k=1}^{\infty}p_{n,k,\gamma}(x) \int_0^{\infty}\left\{(k-1)-(n+2\gamma)t\right\}b_{n,k,\gamma}(t) \left(\frac{nt+\alpha}{n+\beta}-x\right)^m dt\right.\nonumber\\&&\left.+\sum_{k=1}^{\infty}p_{n,k,\gamma}(x) \int_0^{\infty}b_{n,k,\gamma}(t) \left(\frac{nt+\alpha}{n+\beta}-x\right)^m dt \right\}\nonumber\\
&&+\left\{ (n+2\gamma)\sum_{k=1}^{\infty}p_{n,k,\gamma}(x) \int_0^{\infty}tb_{n,k,\gamma}(t) \left(\frac{nt+\alpha}{n+\beta}-x\right)^m dt\right\}:=I_1+I_2 \text{ (say) }.
\end{eqnarray}
To estimate $I_2$ using $t= \frac{n+\beta}{n}\left\{\left(\frac{nt+\alpha}{n+\beta}-x\right)-\left(\frac{\alpha}{n+\beta}-x\right)\right\},$ we have
\begin{eqnarray}
I_2&=&
\frac{(n+2\gamma)(n+\beta)}{n} \left\{\sum_{k=1}^{\infty}p_{n,k,\gamma}(x) \int_0^{\infty}b_{n,k,\gamma}(t) \left(\frac{nt+\alpha}{n+\beta}-x\right)^{m+1}dt\right.\nonumber\\
&&\left.-\left(\frac{\alpha}{n+\beta}-x\right)\sum_{k=1}^{\infty}p_{n,k,\gamma}(x) \int_0^{\infty}b_{n,k,\gamma}(t) \left(\frac{nt+\alpha}{n+\beta}-x\right)^m dt \right\}\nonumber\\
&=& \frac{(n+2\gamma)(n+\beta)}{n}\left\{\sum_{k=1}^{\infty}p_{n,k,\gamma}(x) \int_0^{\infty}b_{n,k,\gamma}(t) \left(\frac{nt+\alpha}{n+\beta}-x\right)^{m+1} dt+ p_{n,0,\gamma}(x)\left(\frac{\alpha}{n+\beta}-x\right)^{m+1}\right.\nonumber\\
&&\left.- \left(\frac{\alpha}{n+\beta}-x\right)\left\{\sum_{k=1}^{\infty}p_{n,k,\gamma}(x) \int_0^{\infty}b_{n,k,\gamma}(t) \left(\frac{nt+\alpha}{n+\beta}-x\right)^{m}dt+ p_{n,0,\gamma}(x)\left(\frac{\alpha}{n+\beta}-x\right)^m\right\}\right\}\nonumber\\
&=&\frac{(n+2\gamma)(n+\beta)}{n} \left\{\mu_{n,m+1,\gamma}(x)-\left(\frac{\alpha}{n+\beta}-x\right) \mu_{n,m,\gamma}(x)\right\}.
\end{eqnarray}
Next, to estimate $I_1$ using the equality, $ \left\{(k-1)-(n+2\gamma)t\right\}b_{n,k,\gamma}(t) = t(1+\gamma t )b_{n,k,\gamma}^{(1)} (t)$, we have
\begin{eqnarray}
I_1&=& \left\{\sum_{k=1}^{\infty}p_{n,k,\gamma}(x) \int_0^{\infty}t b_{n,k,\gamma}^{(1)}(t) \left(\frac{nt+\alpha}{n+\beta}-x\right)^m dt+\sum_{k=1}^{\infty}p_{n,k,\gamma}(x) \int_0^{\infty}b_{n,k,\gamma}(t) \left(\frac{nt+\alpha}{n+\beta}-x\right)^m dt \right\}\nonumber\\
&&+\gamma \sum_{k=1}^{\infty}p_{n,k,\gamma}(x) \int_0^{\infty}t^2 b_{n,k,\gamma}^{(1)}(t) \left(\frac{nt+\alpha}{n+\beta} - x\right)^m dt := J_1+J_2 \text{ (say) }
\end{eqnarray}
Again putting $t= \frac{n+\beta}{n}\left\{\left(\frac{nt+\alpha}{n+\beta}-x\right)-\left(\frac{\alpha}{n+\beta}-x\right)\right\},$ we get
\begin{eqnarray*}
J_1&=& \frac{n+\beta}{n}\left\{\sum_{k=1}^{\infty}p_{n,k,\gamma}(x) \int_0^{\infty} b_{n,k,\gamma}^{(1)}(t) \left(\frac{nt+\alpha}{n+\beta}-x\right)^{m+1} dt\right.\\
&&\left.+\left(\frac{\alpha}{n+\beta}-x \right)\sum_{k=1}^{\infty}p_{n,k,\gamma}(x) \int_0^{\infty}b_{n,k,\gamma}^{(1)}(t) \left(\frac{nt+\alpha}{n+\beta}-x\right)^m dt \right\}\\
&&+\sum_{k=1}^{\infty}p_{n,k,\gamma}(x) \int_0^{\infty}b_{n,k,\gamma}(t) \left(\frac{nt+\alpha}{n+\beta}-x\right)^m dt.
\end{eqnarray*}
Now integrating by parts, we get
\begin{eqnarray}
J_1 &=& -(m+1)\sum_{k=1}^{\infty}p_{n,k,\gamma}(x) \int_0^{\infty} b_{n,k,\gamma}(t) \left(\frac{nt+\alpha}{n+\beta}-x\right)^{m} dt\nonumber\\
&&+m \left(\frac{\alpha}{n+\beta}-x \right)\sum_{k=1}^{\infty}p_{n,k,\gamma}(x) \int_0^{\infty}b_{n,k,\gamma}(t) \left(\frac{nt+\alpha}{n+\beta}-x\right)^{m-1} dt \nonumber\\
&&+\sum_{k=1}^{\infty}p_{n,k,\gamma}(x) \int_0^{\infty}b_{n,k,\gamma}(t) \left(\frac{nt+\alpha}{n+\beta}-x\right)^m dt\nonumber\\
&=& -(m+1)\left\{\sum_{k=1}^{\infty}p_{n,k,\gamma}(x) \int_0^{\infty} b_{n,k,\gamma}(t) \left(\frac{nt+\alpha}{n+\beta}-x\right)^{m} dt+ p_{n,0,\gamma}(x)\left(\frac{\alpha}{n+\beta}-x\right)^m\right\}\nonumber\\
&&+m \left(\frac{\alpha}{n+\beta}-x \right)\left\{\sum_{k=1}^{\infty}p_{n,k,\gamma}(x) \int_0^{\infty}b_{n,k,\gamma}(t) \left(\frac{nt+\alpha}{n+\beta}-x\right)^{m-1} dt \right.+\left. p_{n,0,\gamma}(x)\left(\frac{\alpha}{n+\beta}-x\right)^{m-1}\right\} \nonumber\\
&&+\sum_{k=1}^{\infty}p_{n,k,\gamma}(x) \int_0^{\infty}b_{n,k,\gamma}(t) \left(\frac{nt+\alpha}{n+\beta}-x\right)^m dt+p_{n,0,\gamma}(x)\left(\frac{\alpha}{n+\beta}-x\right)^m\nonumber\\
J_1&=&-m\mu_{n,m,\gamma}(x)+m \left(\frac{\alpha}{n+\beta}-x\right)\mu_{n,m-1,\gamma}(x).
\end{eqnarray}
Proceeding in the similar manner, we obtain the estimate $J_2$ as
\begin{eqnarray}\label{5.eq2.7}
J_2&=& - \frac{\gamma(n+\beta)(m+2)}{n} \mu_{n,m+1,\gamma}(x)+ 2\gamma\frac{(n+\beta)(m+1)}{n}\left(\frac{\alpha}{n+\beta}-x\right)\mu_{n,m,\gamma}(x)\nonumber\\
&&- \frac{m\gamma(n+\beta)}{n}\left(\frac{\alpha}{n+\beta}-x\right)^2\mu_{n,m-1,\gamma}(x).
\end{eqnarray}
Combining (\ref{5.eq2.1})-(\ref{5.eq2.7}), we get\\
$\displaystyle
x(1+\gamma x) \left\{\mu_{n,m,\gamma}^{(1)}(x)+m\mu_{n,m-1,\gamma}(x)\right\}\\
\indent \indent= -m\mu_{n,m,\gamma}(x)+ m \left(\frac{\alpha}{n+\beta}-x\right)\mu_{n,m-1,\gamma}(x)- \frac{\gamma(n+\beta)(m+2)}{n} \mu_{n,m+1,\gamma}(x)\\
\indent \indent\indent \indent + 2\gamma\frac{(n+\beta)(m+1)}{n}\left(\frac{\alpha}{n+\beta}-x\right)\mu_{n,m,\gamma}(x)\\
\indent \indent\indent \indent - \frac{m\gamma(n+\beta)}{n}\left(\frac{\alpha}{n+\beta}-x\right)^2\mu_{n,m-1,\gamma}(x)-n x \mu_{n,m,\gamma}(x)\\
\indent \indent\indent \indent + \frac{(n+2\gamma)(n+\beta)}{n} \left\{\mu_{n,m+1,\gamma}(x)-\left(\frac{\alpha}{n+\beta}-x\right) \mu_{n,m,\gamma}(x)\right\}.
$\\
Hence,
\begin{eqnarray*}
(n-\gamma m)(n+\beta)\mu_{n,m+1,\gamma}(x)&=& n x (1+\gamma x) \left\{\mu_{n,m,\gamma}^{(1)}(x)+m\mu_{n,m-1,\gamma}(x)\right\}\\
&& + \left\{m n+ n^2x- \left(2\gamma m-n\right)\left(\alpha -(n+\beta)x\right)\right\}\mu_{n,m,\gamma}(x)\\
&& + \left\{  m \gamma (n+\beta) \left(\frac{\alpha}{n+\beta}-x \right)^2-m n \left(\frac{\alpha}{n+\beta}-x\right) \right\}\mu_{n,m-1,\gamma}(x).
\end{eqnarray*}
\begin{remark} \label{5.lemma1} \cite{li2005voronovskaja} For $m\in \mathbb{N}\cup \{0\}$, if the m$^{th}$ order moment is defined as
$$ U_{n,m,\gamma} (x) = \sum_{k=0}^{\infty} p_{n,k,\gamma}(x) \left(\frac{k}{n}-x\right)^m,$$
then $U_{n,0,\gamma} (x)= 1$, $U_{n,1,\gamma} (x)= 0$ and
$\displaystyle  nU_{n,m+1,\gamma} (x)= x(1+\gamma x)\left(U_{n,m,\gamma}^{(1)} (x) + m U_{n,m-1,\gamma} (x)\right).$\\
Consequently, for all $x\in [0,\infty)$, we have $U_{n,m,\gamma}
(x)= O(n^{-[(m+1)/2]}).$
\end{remark}
\begin{remark}\label{5.remark2}
It is easily verified from Lemma \ref{5.lemma4} that for each
$x\in [0,\infty)$
\begin{eqnarray*}
B_{n,\gamma}^{\alpha,\beta}(t^m,x)&=&  \frac{n^m\Gamma(n/\gamma +m) \Gamma(n/\gamma -m+1)}{(n+\beta)^m \Gamma(n/\gamma +1)\Gamma(n/\gamma)}x^m \\
&&+ \frac{m n^{m-1} \Gamma(n/\gamma +m-1) \Gamma(n/\gamma -m+1)}{(n+\beta)^m \Gamma(n/\gamma +1)\Gamma(n/\gamma)}\left\{n(m-1)+\alpha (n/\gamma-m+1)\right\}x^{m-1}\\
&&+  \frac{\alpha m(m-1) n^{m-2} \Gamma(n/\gamma +m-2) \Gamma(n/\gamma -m+2)}{(n+\beta)^m \Gamma(n/\gamma +1)\Gamma(n/\gamma)}\\
&&\times\left\{n(m-2)+\frac{\alpha (n/\gamma-m+2)}{2}\right\}x^{m-2}+ O(n^{-2}).
\end{eqnarray*}
\end{remark}
\begin{lemma}\cite{li2005voronovskaja} \label{5.lemma5}
The polynomials $Q_{i,j,r,\gamma}(x)$ exist independent of $n$ and
$k$ such that
$$\{x(1+\gamma x)^r\}D^r[p_{n,k,\gamma}(x)]= \sum_{\substack{2i+j \le r \\ i,j \ge 0}} n^i (k-nx)^j Q_{i,j,r,\gamma}(x)p_{n,k,\gamma}(x), \text{ where }D\equiv \frac{d}{dx}.$$
\end{lemma}
\begin{lemma}\label{5.lemma6}
If $f$ is $r$ times differentiable on $[0,\infty)$, such that
$f^{(r-1)} = O(t^\upsilon), \upsilon>0$ as $t\to \infty$, then for
$r=1,2,3,\ldots$ and $n>\upsilon +\gamma r$ we have
$$ (B_{n,\gamma}^{\alpha,\beta})^{(r)}(f,x)=
\frac{n^r\Gamma(n/\gamma +r) \Gamma(n/\gamma -r+1)}{(n+\beta)^r\Gamma(n/\gamma +1)\Gamma(n/\gamma)}\sum_{k=0}^{\infty}   p_{n+\gamma r,k,\gamma}(x) \int_0^{\infty}b_{n-\gamma r,k+r,\gamma}(t)f^{(r)}\left(\frac{nt+\alpha}{n+\beta}\right)dt.
$$
\end{lemma}
\textbf{Proof:}
First
$$\displaystyle (B_{n,\gamma}^{\alpha,\beta})^{(1)}(f,x)=\sum_{k=1}^{\infty} p_{n,k,\gamma}^{(1)} (x) \int_0^{\infty}b_{n,k,\gamma}(t) f\left(\frac{nt+\alpha}{n+\beta}\right)dt- n(1+\gamma x)^{-n/\gamma-1}f\left(\frac{\alpha}{n+\beta}\right)$$
Now, using the identities
\begin{eqnarray}
 p_{n,k,\gamma}^{(1)} (x)&=& n \left\{p_{n+\gamma,k-1,\gamma}(x) - p_{n+\gamma,k,\gamma}(x)\right\}\label{eq:5.eq2.8}\\
 b_{n,k,\gamma}^{(1)} (x)&=& (n+\gamma) \left\{b_{n+\gamma,k-1,\gamma}(x) - b_{n+\gamma,k,\gamma}(x)\right\},\label{eq:5.eq2.9}
\end{eqnarray}
for $k\ge1$, we have
\begin{eqnarray*}
(B_{n,\gamma}^{\alpha,\beta})^{(1)}(f,x)&=&\sum_{k=1}^{\infty}  n \left\{p_{n+\gamma,k-1,\gamma}(x) - p_{n+\gamma,k,\gamma}(x)\right\} \int_0^{\infty}b_{n,k,\gamma}(t) f\left(\frac{nt+\alpha}{n+\beta}\right)dt\\
&&- n(1+\gamma x)^{-n/\gamma-1}f\left(\frac{\alpha}{n+\beta}\right)\\
&=&  np_{n+\gamma,0,\gamma}(x)\int_0^{\infty}b_{n+\gamma,1,\gamma}( t )f\left(\frac{nt+\alpha}{n+\beta}\right)dt- n(1+\gamma x)^{-n/\gamma-1}f\left(\frac{\alpha}{n+\beta}\right)\\
&&+n\sum_{k=1}^{\infty}   p_{n+\gamma,k,\gamma}(x) \int_0^{\infty}\left\{b_{n,k+1,\gamma}(t)-b_{n,k,\gamma}(t)\right\} f\left(\frac{nt+\alpha}{n+\beta}\right)dt\end{eqnarray*}
\begin{eqnarray*}
(B_{n,\gamma}^{\alpha,\beta})^{(1)}(f,x) &=&  n (1+\gamma x)^{-n/\gamma -1}\int_0^{\infty}(n+\gamma)(1+\gamma t )^{-n/\gamma -2} f\left(\frac{nt+\alpha}{n+\beta}\right)dt\\
&&+n\sum_{k=1}^{\infty}   p_{n+\gamma,k,\gamma}(x) \int_0^{\infty}\left(-\frac{1}{n}b_{n-\gamma,k+1,\gamma}^{(1)}(t)\right)f\left(\frac{nt+\alpha}{n+\beta}\right)dt\\
&&- n(1+\gamma x)^{-n/\gamma-1}f\left(\frac{\alpha}{n+\beta}\right).
\end{eqnarray*}
Integrating by parts, we get\\
$\displaystyle
(B_{n,\gamma}^{\alpha,\beta})^{(1)}(f,x)=n (1+\gamma x)^{-n/\gamma -1}f\left(\frac{\alpha}{n+\beta}\right)
 + \frac{n^2}{n+\beta} (1+\gamma x)^{-n/\gamma -1}\int_0^{\infty}(1+\gamma t )^{-n/\gamma -1} f^{(1)}\left(\frac{nt+\alpha}{n+\beta}\right)dt\\
\indent \indent \indent\indent \indent + \frac{n}{n+\beta}\sum_{k=1}^{\infty}   p_{n+\gamma,k,\gamma}(x) \int_0^{\infty}b_{n-\gamma,k+1,\gamma}(t)f^{(1)}\left(\frac{nt+\alpha}{n+\beta}\right)dt
- n(1+\gamma x)^{-n/\gamma-1}f\left(\frac{\alpha}{n+\beta}\right)\\
(B_{n,\gamma}^{\alpha,\beta})^{(1)}(f,x)
=\frac{n}{n+\beta}\sum_{k=0}^{\infty}   p_{n+\gamma,k,\gamma}(x) \int_0^{\infty}b_{n-\gamma,k+1,\gamma}(t)f^{(1)}\left(\frac{nt+\alpha}{n+\beta}\right)dt. $\\
Thus the result is true for $r=1$. We prove the result by induction method. Suppose that the result is true for $r=i$, then
\begin{eqnarray*}
(B_{n,\gamma}^{\alpha,\beta})^{(i)}(f,x)&=&\frac{n^i\Gamma(n/\gamma +i) \Gamma(n/\gamma -i+1)}{(n+\beta)^i\Gamma(n/\gamma +1)\Gamma(n/\gamma)}\sum_{k=0}^{\infty}   p_{n+\gamma i,k,\gamma}(x) \int_0^{\infty}b_{n-\gamma i,k+i,\gamma}(t)f^{(i)}\left(\frac{nt+\alpha}{n+\beta}\right)dt.
\end{eqnarray*}
Thus using the identities \eqref{eq:5.eq2.8} and \eqref{eq:5.eq2.9}, we have\\
$\displaystyle
(B_{n,\gamma}^{\alpha,\beta})^{(i+1)}(f,x)\\
\indent \indent=\frac{n^i\Gamma(n/\gamma +i) \Gamma(n/\gamma -i+1)}{(n+\beta)^i\Gamma(n/\gamma +1)\Gamma(n/\gamma)}\left\{\sum_{k=1}^{\infty} (n/\gamma + i)  \left\{p_{n+\gamma (i+1),k-1,\gamma}(x)- p_{n+\gamma (i+1),k,\gamma}(x)\right\}\right. \\
\indent \indent \indent \left.\times \int_0^{\infty}b_{n-\gamma i,k+i,\gamma}(t)f^{(i)}\left(\frac{nt+\alpha}{n+\beta}\right)dt\right.\\
\left. \indent \indent \indent - \left(n/\gamma+ i\right) ( 1+\gamma x)^{-n/\gamma -i-1} \int_0^{\infty} b_{n-\gamma i,i,\gamma}(t) f^{(i)}\left(\frac{nt+\alpha}{n+\beta}\right)\right\}\\
\indent \indent =\frac{n^i\Gamma(n/\gamma +i+1) \Gamma(n/\gamma -i+1)}{(n+\beta)^i\Gamma(n/\gamma +1)\Gamma(n/\gamma)}p_{n+\gamma (i+1),0,\gamma}(x)\int_0^{\infty}b_{n-\gamma i,1+i,\gamma}(t)f^{(i)}\left(\frac{nt+\alpha}{n+\beta}\right)dt \\
\indent \indent \indent  -\frac{n^i\Gamma(n/\gamma +i+1) \Gamma(n/\gamma -i+1)}{(n+\beta)^i\Gamma(n/\gamma +1)\Gamma(n/\gamma)}p_{n+\gamma (i+1),0,\gamma}(x)\int_0^{\infty}b_{n-\gamma i,i,\gamma}(t)f^{(i)}\left(\frac{nt+\alpha}{n+\beta}\right)dt \\
\indent \indent \indent  + \frac{n^i\Gamma(n/\gamma +i+1) \Gamma(n/\gamma -i+1)}{(n+\beta)^i\Gamma(n/\gamma +1)\Gamma(n/\gamma)}\left\{\sum_{k=1}^{\infty} p_{n+\gamma (i+1),k,\gamma}(x)\right.\\
\indent \indent \indent \left.\times \int_0^{\infty}\left\{b_{n-\gamma i,k+i+1,\gamma}(t)-b_{n-\gamma i,k+i,\gamma}(t)\right\}f^{(i)}\left(\frac{nt+\alpha}{n+\beta}\right)dt\right\}\\
\indent \indent =\frac{n^i\Gamma(n/\gamma +i+1) \Gamma(n/\gamma -i+1)}{(n+\beta)^i\Gamma(n/\gamma +1)\Gamma(n/\gamma)}p_{n+\gamma (i+1),0,\gamma}(x)\\
\indent \indent \indent \times \int_0^{\infty}\left(-\frac{1}{n/\gamma -i}b_{n-\gamma (i-1),1+i,\gamma}^{(1)}(t)\right)f^{(i)}\left(\frac{nt+\alpha}{n+\beta}\right)dt \\
\indent \indent \indent + \frac{n^i\Gamma(n/\gamma +i+1) \Gamma(n/\gamma -i+1)}{(n+\beta)^i\Gamma(n/\gamma +1)\Gamma(n/\gamma)}\sum_{k=1}^{\infty} p_{n+\gamma (i+1),k,\gamma}(x)\\
\indent \indent \indent \times \int_0^{\infty}\left(-\frac{1}{n/\gamma -i}b_{n-\gamma (i-1),k+i+1,\gamma}^{(1)}(t)\right)f^{(i)}\left(\frac{nt+\alpha}{n+\beta}\right)dt.
$\\
Integrating by parts, we obtain
\begin{eqnarray*}
(B_{n,\gamma}^{\alpha,\beta})^{(i+1)}(f,x)&=&\frac{n^{i+1}\Gamma(n/\gamma +i+1) \Gamma(n/\gamma -i+1)}{(n+\beta)^{i+1}\Gamma(n/\gamma +1)\Gamma(n/\gamma)}\sum_{k=0}^{\infty} p_{n+\gamma (i+1),k,\gamma}(x)\\&&\times\int_0^{\infty}b_{n-\gamma (i-1),k+i+1,\gamma}(t)f^{(i+1)}\left(\frac{nt+\alpha}{n+\beta}\right)dt.
\end{eqnarray*}
This completes the proof of Lemma \ref{5.lemma6}.
\section{Direct Theorems}
This section deals with the direct results, we establish here pointwise approximation, asymptotic formula and error estimation in simultaneous approximation.\\
\indent  We denote $C_{\mu} [0,\infty)= \{ f\in C[0,\infty):
|f(t)|\leq Mt^{\mu}$ for some $ M>0, \mu>0\},$ and the norm
$\|\cdot\|_{\mu}$ on the class $C_{\mu} [0,\infty)$ is defined as
$\displaystyle \|f\|_{\mu}= \sup_{0\le t < \infty}
|f(t)|t^{-\mu}.$ It can be easily verified that the operators
$B_{n,\gamma}^{\alpha,\beta}(f,x)$ are well defined for $f\in
C_{\mu}[0,\infty)$.
\begin{theorem}\label{5.thm1}
Let $f\in C_{\mu} [0,\infty)$ and $f^{(r)}$ exists at a point $x\in (0,\infty)$. Then we have
$$\lim_{n \to \infty} \left(B_{n,\gamma}^{\alpha,\beta}\right)^{(r)} (f,x) = f^{(r)}(x).$$
\end{theorem}
\textbf{Proof:}
By Taylor's expansion of $f$, we have
$$f(t) = \sum_{i=0}^r \frac{f^{(i)}(x)}{i!} (t-x)^i + \epsilon(t,x) (t-x)^r,$$
where $\epsilon(t,x) \to 0$ as $t\to x$. Hence,
\begin{eqnarray*}
\left(B_{n,\gamma}^{\alpha,\beta}\right)^{(r)}(f,x) &=& \sum_{i=0}^r \frac{f^{(i)}(x)}{i!}\left(B_{n,\gamma}^{\alpha,\beta}\right)^{(r)}((t-x)^i,x)+ \left(B_{n,\gamma}^{\alpha,\beta}\right)^{(r)}(\epsilon(t,x)(t-x)^r,x)\\
&:=& R_1+ R_2.
\end{eqnarray*}
First to estimate $R_1$, using binomial expansion of $\left(\frac{nt+\alpha}{n+\beta}-x\right)^i$ and Remark \ref{5.remark2}, we have
\begin{eqnarray*}
R_1 &=& \sum_{i=0}^r \frac{f^{(i)}(x)}{i!} \sum_{j=0}^i \binom{i}{j} (-x)^{i-j} \left(B_{n,\gamma}^{\alpha,\beta}\right)^{(r)}(t^j,x)\\
&=& \frac{f^{(r)}(x)}{r!}\left\{ \frac{n^r \Gamma(n/\gamma +r)\Gamma(n/\gamma -r+1)}{(n+\beta)^r\Gamma(n/\gamma +1)\Gamma(n/\gamma)}r! \right\}\\
&=&  f^{(r)}(x)\left\{ \frac{n^r \Gamma(n/\gamma +r)\Gamma(n/\gamma -r+1)}{(n+\beta)^r\Gamma(n/\gamma +1)\Gamma(n/\gamma)}\right\}\to f^{(r)}(x)\text{ as } n\to \infty.
\end{eqnarray*}
Next applying Lemma \ref{5.lemma5}, we obtain
$$R_2 = \int_0^{\infty} W_{n,\gamma}^{(r)}(t,x) \epsilon(t,x) \left(\frac{nt+\alpha}{n+\beta} -x\right)^r dt,$$
\begin{eqnarray*}
|R_2| &\le& \sum_{\substack{2i+j \le r\\  i,j \ge 0}} n^i \frac{|Q_{i,j,r,\gamma}(x)|}{\{x(1+\gamma x)\}^r} \sum_{k=1}^{\infty} |k-n x|^j p_{n,k,\gamma}(x) \int_0^{\infty} b_{n,k,\gamma}(t) |\epsilon(t,x)|\bigg|\frac{nt+\alpha}{n+\beta}-x\bigg|^r dt\\
 &&+ \frac{\Gamma(n/\gamma +r+2)}{\Gamma(n/\gamma)}(1+\gamma x)^{-n/\gamma-r} |\epsilon(0,x)| \bigg|\frac{\alpha}{n+\beta}-x\bigg|^r.
\end{eqnarray*}
The second term in the above expression tends to zero as $n\to \infty$.
Since $\epsilon(t,x) \to 0$ as $t\to x $ for given $\varepsilon>0,$ there exists a $\delta\in (0,1)$
 such that $|\epsilon(t,x)| < \varepsilon$ whenever $0<|t-x| < \delta$.
 If $\tau> \max\{ \mu, r\}$, where $\tau$ is any integer, then we can find a constant $M_3>0$, such that $|\epsilon(t,x)\big(\frac{nt+\alpha}{n+\beta}-x\big)^r|\le M_3\big|\frac{nt+\alpha}{n+\beta}-x\big|^\tau$ for $|t-x| \ge \delta$. Therefore,\\
$\displaystyle
|R_2|  \le  M_3 \sum_{\substack{2i+j \le r \\ i,j \ge 0}} n^i\sum_{k=0}^{\infty}  |k-n x|^j p_{n,k,\gamma}(x)\left\{\varepsilon \int_{|t-x|<\delta} b_{n,k,\gamma}(x) \bigg|\frac{nt+\alpha}{n+\beta}-x\bigg|^r dt\right.\\
\indent \indent \indent \left.+ \int_{|t-x| \ge \delta} b_{n,k,\gamma}(t) \bigg|\frac{nt+\alpha}{n+\beta}-x\bigg|^\tau dt\right\}=: R_3 + R_4.
$\\
Applying the Cauchy-Schwarz inequality for integration and summation, respectively, we obtain\\
$\displaystyle  |R_3|\leq  \varepsilon M_3 \sum_{\substack{2i+j \le r \\ i,j \ge 0}} n^i\left\{\sum_{k=1}^{\infty}  (k-n x)^{2j} p_{n,k,\gamma}(x) \right\}^{1/2} \left\{ \sum_{k=1}^{\infty} p_{n,k,\gamma}(x)\int_{0}^{\infty} b_{n,k,\gamma}(t) \left(\frac{nt+\alpha}{n+\beta}-x\right)^{2r} dt\right\}^{1/2}.$\\
Using Remark \ref{5.lemma1} and Lemma \ref{5.lemma4}, we get
$R_3\leq  \varepsilon O(n^{r/2}) O(n^{-r/2}) = \varepsilon \cdot O(1).$\\
Again using the Cauchy-Schwarz inequality and Lemma \ref{5.lemma4}, we get
\begin{eqnarray*}
|R_4| &\leq& M_4 \sum_{\substack{2i+j \le r \\ i,j \ge 0}} n^i\sum_{k=1}^{\infty}  |k-n x|^{j} p_{n,k,\gamma}(x)\int_{|t-x|\ge \delta} b_{n,k,\gamma}(t) \bigg|\frac{nt+\alpha}{n+\beta}-x\bigg|^{\tau} dt\\
&\leq& M_4 \sum_{\substack{2i+j \le r \\ i,j \ge 0}} n^i\sum_{k=1}^{\infty}  |k-n x|^{j} p_{n,k,\gamma}(x)\left\{\int_{|t-x|\ge \delta} b_{n,k,\gamma}(t) dt\right\}^{1/2}\\
&&\times \left\{\int_{|t-x|\ge \delta} b_{n,k,\gamma}(t) \bigg(\frac{nt+\alpha}{n+\beta}-x\bigg)^{2\tau} dt\right\}^{1/2}\\
&\leq& M_4 \sum_{\substack{2i+j \le r\\ i,j \ge 0}} n^i \left\{\sum_{k=1}^{\infty}  (k-n x)^{2j} p_{n,k,\gamma}(x)\right\}^{1/2}\left\{\sum_{k=1}^{\infty} p_{n,k,\gamma}(x)\int_{0}^{\infty} b_{n,k,\gamma}(t) \bigg(\frac{nt+\alpha}{n+\beta}-x\bigg)^{2\tau} dt\right\}^{1/2}\\
&=& \sum_{\substack{2i+j \le r \\ i,j \ge 0}} n^iO(n^{j/2})O(n^{-\tau/2})= O(n^{(r-\tau)/2})= o(1).
\end{eqnarray*}
Collecting the estimation of $R_1-R_4$, we get the required result.
\begin{theorem}\label{5.thm2}
Let $f\in C_{\mu}[0,\infty)$. If $f^{(r+2)}$ exists at a point $x\in (0,\infty)$, then
 \begin{eqnarray}\label{5.eq3.1}\lim_{n\to \infty} n \left\{\left(B_{n,\gamma}^{\alpha,\beta}\right)^{(r)}(f,x) - f^{(r)}(x)\right\}&=& r(\gamma (r-1) -\beta)f^{(r)}(x) + \left\{r\gamma(1+2x)+ \alpha-\beta x\right\}f^{(r+1)}(x)\nonumber\\
 &&+ x(1+\gamma x ) f^{(r+2)}(x).\end{eqnarray}
\end{theorem}
\textbf{Proof:}
Using Taylor's expansion of $f$, we have
$$\displaystyle f(t) = \sum_{i=0}^{r+2} \frac{f^{(i)}(x)}{i!}(t-x)^i + \epsilon(t,x) (t-x)^{r+2},$$
where $\epsilon(t,x) \to 0$ as $t\to x $ and $\epsilon(t,x) = O((t-x)^{\mu}), t \to \infty$ for $\mu>0$.\\
Applying Lemma \ref{5.lemma4}, we have\\
$\displaystyle  n \left\{\left(B_{n,\gamma}^{\alpha,\beta}\right)^{(r)}(f,x) - f^{(r)}(x)\right\} = n \left\{\sum_{i=0}^{r+2} \frac{f^{(i)}(x)}{i!}\left(B_{n,\gamma}^{\alpha,\beta}\right)^{(r)}((t-x)^i,x) -f^{(r)}(x)\right\}\\
\indent\indent \indent \indent\indent \indent \indent\indent \indent + n\left\{\left(B_{n,\gamma}^{\alpha,\beta}\right)^{(r)}(\epsilon(t,x)(t-x)^{r+2},x)\right\}:= E_1+E_2
$\\
First, we have
\begin{eqnarray*}
E_1&=& n \sum_{i=0}^{r+2} \frac{f^{(i)}(x)}{i!} \sum_{j=0}^i \binom{i}{j} (-x)^{i-j}\left(B_{n,\gamma}^{\alpha,\beta}\right)^{(r)}(t^j,x) -nf^{(r)}(x)\\
&=& \frac{f^{(r)}(x)}{r!}n\left\{\left(B_{n,\gamma}^{\alpha,\beta}\right)^{(r)}(t^r,x)-r!\right\}\\
&&+ \frac{f^{(r+1)}(x)}{(r+1)!}n\left\{(r+1)(-x) \left(B_{n,\gamma}^{\alpha,\beta}\right)^{(r)}(t^r,x)+\left(B_{n,\gamma}^{\alpha,\beta}\right)^{(r)}(t^{r+1},x)\right\}\\
&& + \frac{f^{(r+2)}(x)}{(r+2)!}n\left\{\frac{(r+2)(r+1)}{2}x^2 \left(B_{n,\gamma}^{\alpha,\beta}\right)^{(r)}(t^r,x)+(r+2)(-x) \left(B_{n,\gamma}^{\alpha,\beta}\right)^{(r)}(t^{r+1},x)\right.\\
&&\left.+\left(B_{n,\gamma}^{\alpha,\beta}\right)^{(r)}(t^{r+2},x)\right\}\\
&=& f^{(r)}(x)n\left\{\frac{n^r\Gamma(n/\gamma+r)\Gamma(n/\gamma -r +1)}{(n+\beta)^r\Gamma(n/\gamma+1)\Gamma(n/\gamma)}-1\right\}\\
&&+ \frac{f^{(r+1)}(x)}{(r+1)!}n\left\{(r+1)(-x)\frac{n^r\Gamma(n/\gamma+r)\Gamma(n/\gamma -r +1)}{(n+\beta)^r\Gamma(n/\gamma+1)\Gamma(n/\gamma)}r!\right.\\
&&\left.+ \frac{n^{r+1}\Gamma(n/\gamma+r+1)\Gamma(n/\gamma -r)}{(n+\beta)^{r+1}\Gamma(n/\gamma+1)\Gamma(n/\gamma)} (r+1)! x\right.\\
&&\left. + \frac{(r+1)n^r\Gamma(n/\gamma+r)\Gamma(n/\gamma -r)}{(n+\beta)^{r+1}\Gamma(n/\gamma+1)\Gamma(n/\gamma)}\left\{nr+\alpha(n/\gamma-r)\right\}r!\right\}\\
&&+\frac{f^{(r+2)}(x)}{(r+2)!}n\left( \frac{(r+1)(r+2)}{2} x^2 \frac{n^r\Gamma(n/\gamma+r)\Gamma(n/\gamma -r +1)}{(n+\beta)^r\Gamma(n/\gamma+1)\Gamma(n/\gamma)}r! \right.\\
&&\left.-x(r+2)\left\{\frac{n^{r+1}\Gamma(n/\gamma+r+1)\Gamma(n/\gamma -r)}{(n+\beta)^{r+1}\Gamma(n/\gamma+1)\Gamma(n/\gamma)}(r+1)!x \right.\right.\\
&&\left.\left.+ \frac{(r+1)n^r\Gamma(n/\gamma+r)\Gamma(n/\gamma -r)}{(n+\beta)^{r+1}\Gamma(n/\gamma+1)\Gamma(n/\gamma)}\left\{nr+\alpha(n/\gamma-r)\right\}r!\right\}\right.\\
&&\left.+\frac{n^{r+2}\Gamma(n/\gamma+r+2)\Gamma(n/\gamma -r -1)}{(n+\beta)^{r+2}\Gamma(n/\gamma+1)\Gamma(n/\gamma)}\frac{(r+2)!}{2}x^2\right.\\
&&\left.+\frac{(r+2)n^{r+1}\Gamma(n/\gamma+r+1)\Gamma(n/\gamma -r - 1)}{(n+\beta)^{r+2}\Gamma(n/\gamma+1)\Gamma(n/\gamma)}\left\{n(r+1)+\alpha(n/\gamma-r-1)\right\}(r+1)!x\right.\\
&&\left.+ \frac{\alpha (r+1)(r+2)n^r\Gamma(n/\gamma+r)\Gamma(n/\gamma -r)}{(n+\beta)^r\Gamma(n/\gamma+1)\Gamma(n/\gamma)}\left\{nr+\frac{\alpha(n/\gamma-r)}{2}\right\}r!\right).
\end{eqnarray*}
Now, the coefficients of $f^{(r)}(x)$, $f^{(r+1)}(x)$ and
$f^{(r+2)}(x)$ in the above expression are tend to
$r(\gamma(r-1)-\beta)$, $r\gamma(1+2x)+\alpha-\beta x$ and
$x(1+\gamma x)$ respectively, which follows by using induction
hypothesis on $r$ and taking the limit as $n\to \infty$. Hence in
order to prove \eqref{5.eq3.1}, it is sufficient to show that $E_2
\to 0$ as $n\to \infty$, which follows along the lines of the
proof of Theorem \ref{5.thm1} and by using Remark \ref{5.lemma1},
Lemma \ref{5.lemma4} and Lemma \ref{5.lemma5}.
\begin{remark}
   Particular case $\alpha=\beta=0$ was discussed in (\cite{Vijaygupta2009}, Th. 4.1),
   which says that the coefficient of $f^{(r+1)}(x)$ converges to $r(1+2\gamma x)$
   but it converges to $r\gamma(1+2x)$ and we get this by putting $\alpha=\beta=0$ in above Theorem.
\end{remark}
\begin{definition}
The m$^{th}$ order modulus of continuity $\omega_m(f,\delta,[a,b])$ for a function continuous on [a,b] is defined by
$$\omega_m(f,\delta,[a,b]) = \sup \{|\Delta_h^mf(x)|: |h| \leq \delta; x, x+h \in [a,b]\}.$$
For $m=1,$ $\omega_m(f,\delta)$ is usual modulus of continuity.
\end{definition}
\begin{theorem}
Let $f\in C_{\mu}[0,\infty)$ for some $\mu>0$ and $0<a<a_1<b_1<b<\infty$. Then for $n$ sufficiently large, we have
$$\|(B_{n,\gamma}^{\alpha,\beta})^{(r)}(f,\cdot)- f^{(r)}\|_{C[a_1,b_1]} \leq M_1\omega_2(f^{(r)},n^{-1/2},[a_1,b_1])+M_2n^{-1}\|f\|_{\mu},$$
where $M_1= M_1(r)$ and $M_2= M_2(r,f)$.
\end{theorem}
\textbf{Proof:} Let us assume that $0 < a < a_1 < b_1 < b
<\infty$. For sufficiently small $\eta>0$, we define the function
$f_{\eta,2}$ corresponding to $f\in C_{\mu}[a,b]$ and $t\in
[a_1,b_1]$ as follows:
$$ f_{\eta,2}(t) = \eta^{-2} \int_{-\eta/2}^{\eta/2}  \int_{-\eta/2}^{\eta/2} \left( f(t) - \Delta^2_{h}f(t) \right)dt_1dt_2,$$
where $h=(t_1+t_2)/2$ and $\Delta_{h}^2$ is the second order
forward difference operator with step length $h.$ For $f\in
C[a,b]$, the functions $f_{\eta,2}$ are known as the Steklov mean
of order 2, which satisfy the following properties
\cite{freudpopv1969}:
\begin{enumerate}
\item [(a)] $f_{\eta,2}$ has continuous derivatives up to order 2 over $[a_1,b_1]$;
\item [(b)]$\|f_{\eta,2}^{(r)}\|_{C[a_1,b_1]} \le \hat{M_1} \eta^{-r} \omega_2(f,\eta,[a,b]), r=1,2$;
\item [(c)]$\|f - f_{\eta,2}\|_{C[a_1,b_1]} \le \hat{M_2}\omega_2(f,\eta,[a,b])$;
\item [(d)]$\|f_{\eta,2}\|_{C[a_1,b_1]} \le \hat{M_3}\|f\|_{\mu}$,
\end{enumerate}
where $\hat{M_i}$, $i=1,2,3$ are certain constants which are different in each occurrence and are independent of $f$ and $\eta$. \\
We can write by linearity properties of  $B_{n,\gamma}^{\alpha,\beta}$,
\begin{eqnarray*}
\|(B_{n,\gamma}^{\alpha,\beta})^{(r)}(f,\cdot)- f^{(r)}\|_{C[a_1,b_1]}&\leq&  \|(B_{n,\gamma}^{\alpha,\beta})^{(r)}(f-f_{\eta,2},\cdot)\|_{C[a_1,b_1]}\\
&& +\|(B_{n,\gamma}^{\alpha,\beta})^{(r)}f_{\eta,2},\cdot)-f_{\eta,2}^{(r)}\|_{C[a_1,b_1]}+ \|f^{(r)}-f_{\eta,2}^{(r)}\|_{C[a_1,b_1]}\\
&:=& P_1+ P_2+ P_3.
\end{eqnarray*}
Since $f_{\eta,2}^{(r)} = (f^{(r)})_{\eta,2}(t),$ by property (c)
of the function $f_{\eta,2}$, we get
$$P_3 \le \hat{M_4} \omega_2(f^{(r)},\eta, [a,b]).$$
Next, on an application of Theorem \ref{5.thm2}, it follows that
$$P_2 \leq \hat{M_5} n^{-1} \sum_{i=r}^{r+2} \|f_{\eta,2}^{(i)}\|_{C[a,b]}.$$
Using the interpolation property due to Goldberg and Meir
\cite{Goldberg1971}, for each $j= r, r+1, r+2$, it follows that
$$\|f_{\eta,2}^{(i)}\|_{C[a,b]}\leq \hat{M_6} \{ \|f_{\eta,2}\|_{C[a,b]}  + \|f_{\eta,2}^{(r+2)}\|_{C[a,b]}\}.$$
Therefore, by applying properties (c) and (d) of the function $f_{\eta,2},$ we obtain
$$P_2 \leq \hat{M_7}\cdot n^{-1} \{ \|f\|_{\mu}+ \delta^{-2} \omega_2(f^{(r)},\mu,[a,b])\}.$$
Finally we shall estimate $P_1$, choosing $a^{*}$, $b^{*}$ satisfying the conditions $0< a< a^{*}< a_1< b_1 < b^{*} < b < \infty$. Suppose $\hbar(t)$ denotes the characteristic function of the interval $[a^{*},b^{*}]$. Then
\begin{eqnarray*}
P_1 &\leq& \|(B_{n,\gamma}^{\alpha,\beta})^{(r)}(\hbar(t) \left(f(t)-f_{\eta,2}(t)\right),\cdot)\|_{C[a_1,b_1]}+\|(B_{n,\gamma}^{\alpha,\beta})^{(r)}((1-\hbar(t)) \left(f(t)-f_{\eta,2}(t)\right),\cdot)\|_{C[a_1,b_1]}\\&:=& P_4 + P_5.
\end{eqnarray*}
By Lemma \ref{5.lemma6}, we have \\
$\displaystyle
(B_{n,\gamma}^{\alpha,\beta})^{(r)}(\hbar(t) \left(f(t)-f_{\eta,2}(t)\right),x)\\
\indent \indent = \frac{n^{r}\Gamma(n/\gamma+r)\Gamma(n/\gamma -r +1)}{(n+\beta)^{r}\Gamma(n/\gamma+1)\Gamma(n/\gamma)}\sum_{k=0}^{\infty} p_{n+\gamma r,k,\gamma}(x) \\
\indent \indent \indent \times \int_0^{\infty} b_{n-\gamma r,k+r,\gamma}(t) \hbar(t) \left(f^{(r)}\left(\frac{nt+\alpha}{n+\beta}\right)- f^{(r)}_{\eta,2} \left(\frac{nt+\alpha}{n+\beta}\right)\right)dt
$\\
Hence,
$$\|(B_{n,\gamma}^{\alpha,\beta})^{r}(\hbar(t) \left(f(t)-f_{\eta,2}(t)\right),\cdot)\|_{C[a_1,b_1]}\leq \hat{M_8} \|f^{(r)} - f^{(r)}_{\eta,2}\|_{C[a^{*},b^{*}]}.$$
Now, for $x\in [a_1,b_1]$ and $t\in [0,\infty) \setminus [a^{*},b^{*}]$, we choose a $\delta>0$ satisfying $\bigg|\displaystyle \frac{nt+\alpha}{n+\beta}-x\bigg|\ge \delta.$\\
Therefore, by Lemma \ref{5.lemma5} and the Cauchy-Schwarz
inequality, we have
$$I\equiv (B_{n,\gamma}^{\alpha,\beta})^{(r)}((1-\hbar(t)) \left(f(t)-f_{\eta,2}(t)\right),x)$$
and
\begin{eqnarray*}
|I| &\leq& \sum_{\substack{2i+j \leq r \\i,j \geq 0}}n^i \frac{|Q_{i,j,r,\gamma}(x)|}{\{x(1+\gamma x)\}^r} \sum_{k=1}^{\infty} p_{n,k,\gamma}(x) |k-n x|^j\\
&& \times \int_0^{\infty} b_{n,k,\gamma}(t) \left(1-\hbar(t)\right) \bigg|f\left(\frac{nt+\alpha}{n+\beta}\right) - f_{\eta,2}\left(\frac{nt+\alpha}{n+\beta}\right)\bigg| dt\\
 &&+ \frac{\Gamma(n/\gamma+r)}{\Gamma(n/\gamma)}(1+\gamma x)^{-n/\gamma -r} (1-\hbar(0))\bigg|f\left(\frac{\alpha}{n+\beta}\right)-f_{\eta,2}\left(\frac{\alpha}{n+\beta}\right)\bigg|.
\end{eqnarray*}
For sufficiently large $n$, the second term tends to zero. Thus,
\begin{eqnarray*}
  |I| &\leq& \hat{M_9} \|f\|_{\mu} \sum_{\substack{2i+j \leq r\\ i,j \geq 0}}n^i \sum_{k=1}^{\infty} p_{n,k,\gamma}(x) |k-n x|^j
 \int_{|t-x|\geq \delta} b_{n,k,\gamma}(t) dt\\
&\leq & \hat{M_9} \|f\|_{\mu}\delta^{-2m}  \sum_{\substack{2i+j \leq r \\i,j \geq 0}}n^i \sum_{k=1}^{\infty} p_{n,k,\gamma}(x) |k-n x|^j\\
&&\times \left(\int^{\infty}_0 b_{n,k,\gamma}(t) dt\right)^{1/2}\left(\int^{\infty}_0 b_{n,k,\gamma}(t)\left(\frac{nt+\alpha}{n+\beta}-x\right)^{4m} dt\right)^{1/2}\\
&\leq & \hat{M_9} \|f\|_{\mu}\delta^{-2m}  \sum_{\substack{2i+j \leq r\\ i,j \geq 0}} n^i
 \left(\sum_{k=1}^{\infty}  p_{n,k,\gamma}(x) |k-n x|^{2j}\right)^{1/2}\\
 &&\times\left(\sum_{k=1}^{\infty}  p_{n,k,\gamma}(x)\int^{\infty}_0 b_{n,k,\gamma}(t)\left(\frac{nt+\alpha}{n+\beta}-x\right)^{4m} dt\right)^{1/2}
 \end{eqnarray*}
 Hence, by using Remark \ref{5.lemma1} and Lemma \ref{5.lemma4}, we have
 $$|I|\leq \hat{M}_{10} \|f\|_{\mu} \delta^{-2m} O\left(n^{(i+(j/2) -m)}\right) \leq \hat{M}_{11} n^{-q} \|f\|_{\mu},$$
 where $q= m-(r/2).$ Now choosing $m>0$ satisfying $q\geq 1$, we obtain $I\leq \hat{M}_{11}n^{-1}\|f\|_{\mu}.$ Therefore, by property (c) of the function $f_{\eta,2}(t),$ we get
 \begin{eqnarray*}
 P_1 &\leq& \hat{M_8}\|f^{(r)}- f^{(r)}_{\eta,2} \|_{C[a^{*},b^{*}]} + \hat{M}_{11}n^{-1} \|f\|_{\mu}\\
 &\leq & \hat{M}_{12} \omega_2 (f^{(r)},\eta,[a,b]) + \hat{M}_{11} n^{-1} \|f\|_{\mu}.
 \end{eqnarray*}
 Choosing $\eta = n^{-1/2}$, the theorem follows.
 \begin{remark}
 In the last decade the applications of $q$-calculus in approximation theory is one of the main area of research.
 In 2008, Gupta, \cite{Gupta2008some} is introduced $q$-Durrmeyer operators whose approximation properties were studied in \cite{AralGupta2010}.
 More work in this direction can be seen in \cite{phillips1996bernstein,ostrovska2006lupas,mishra2014generalized}.\\
 A Durrmeyer type $q$-analogue of the $B_{n,\gamma}^{\alpha,\beta}(f,x)$ introduce as follows:
 \begin{equation}\label{qoperator1}
B_{n,\gamma, q}^{\alpha,\beta}(f,x)= \sum_{k=1}^{\infty} p_{n,k,\gamma}^q(x) \int_0^{\infty/A}q^{-k} b_{n,k,\gamma}^q(t) f\left(\frac{[n]_qt+\alpha}{[n]_q+\beta}\right)d_qt+ p_{n,0,\gamma}^q(x) f\left(\frac{\alpha}{[n]_q+\beta}\right),
\end{equation}
where $$p_{n,k,\gamma}^q(x)= q^{\frac{k^2}{2}}\frac{\Gamma_q(n/\gamma + k)}{\Gamma_q(k+1)\Gamma_q(n/\gamma)}\cdot \frac{(q\gamma x)^k}{(1+q\gamma x)^{(n/\gamma)+k}_q}$$
$$b_{n,k,\gamma}^q(x)= \gamma q^{\frac{k^2}{2}}\frac{\Gamma_q(n/\gamma + k+1)}{\Gamma_q(k)\Gamma_q(n/\gamma+1)}\cdot \frac{(\gamma t)^{k-1}}{(1+\gamma t)^{(n/\gamma)+k+1}_q}$$
and
$$\int_0^{\infty/A} f(x) d_q x = (1-q) \sum_{n=-\infty}^{\infty} f\left(\frac{q^n}{A}\right)\frac{q^n}{A},~~~ A>0.$$
  Notations used in (\ref{qoperator1}) can be found in \cite{kac2002quantum}. For the operators (\ref{qoperator1}),
  one can study  their approximation properties based on $q$-integers.
 \end{remark}
 \section*{Competing Interests}
The authors declare that there is no conflict of interests regarding the publication of this research article.
\section*{Acknowledgements}
The authors would like to express their deep gratitude to the
anonymous learned referee(s) and the editor for their valuable
suggestions and constructive comments, which resulted in the
subsequent improvement of this research article.

\end{document}